# Designing videos for undergraduate mathematics

**A Bouzelmate[1] and B Rittaud[2]**

[1] LaR2A laboratory, Faculty of Sciences, Abdelmalek Essaâdi University, Tétouan, Morocco
[2] Laboratoire Analyse, Géométrie et Applications (CNRS UMR 7539), Institut Galilée, Université Sorbonne Paris Nord, Villetaneuse, France

E-mail: rittaud@math.univ-paris13.fr (correspondence author)

**Abstract.** This article presents some technical and pedagogical considerations about a series of videos made by the authors to teach elementary theory of differential equations. Intended to be a pedagogical support for undergraduate students of scientific curricula, these videos provide rigorous mathematical contents in a form that takes into account the specificities of the tool. Videos for teaching mathematics should be seen as a specific way for providing rigorous scientific content to learners, and should not be reduced to oral manuals or standard classroom lessons.

**Research Contribution:** Arij Bouzelmate and Benoît Rittaud are searcher and teacher in mathematics. Arij Bouzelmate is a specialist in partial differential equations and director of LaR2A. Benoît Rittaud is a specialist in discrete mathematics and is involved in popularization of mathematics.

## 1.   Introduction: videos for teaching

Even before the covid-19 pandemics, videos were regarded as a possible pedagogical tool, already widely used in several areas. Now that it seems that the unfolding of digital technologies will only grow bigger and bigger in the foreseeable future, it is reasonable to assert that videos will become more and more used for pedagogical purposes.

Well-designed videos for teaching can be a highly powerful ally in many respects. Their cost is very low when compared to the number of potential watchers. Also, they can be very inclusive by the use of free platforms of access. It can help the promotion of international cooperation between universities, not only of same language-speaking countries since translation of videos could be done at reasonable cost. It is a modern tool for acquisition of knowledge, suited for the youth. It is not intended to replace classroom lessons or written manuals but to complement them. Its intrinsic limitations (the most important being that it makes the student or the schoolchild especially passive) implies that it should not even remotely be regarded as a first step that would make human contacts or exchanges between teachers and learners disappear.

Two common confusions should be avoided when considering videos for pedagogical purposes. The first one is to reduce them as some kind of standard classroom lesson. Sure, practical considerations can easily lead to set up videos by simply putting a camera in the classroom during a lesson. The point is that it is very different, from the learner's point of view, to listen to a teacher in the classroom and to watch him on a video. Especially nowadays where there are so many videos of high quality on so many different subjects, it is very unlikely that a student will willingly watch a video in which so many defects will be heard, as there are inevitably in any oral presentation: hesitations, back an forth explanations because of omissions, slowness, unexpected pauses and timeouts, etc.

The second confusion consists in designing a video in the same way as we would for a manual. Such a confusion typically leads to make the text appear on the screen as it is oralized, either by a presenter or a voice-over. This kind of video is of low interest compared to a classical manual, since, in the latter, the eyes of the reader can catch much more efficiently and more quickly any relevant information than a video would allow.

It is therefore our opinion that videos should avoid to be designed so as to mimic traditional classroom lesson as well as manuals. These are a specific tool, with specific rules and specific interset for learning. They will not replace other available tools, especially because of their intrinsic limitation of not allowing direct interactivity with students, but can still be useful for teachers and learners. We present here some ideas that guided us in the conception of our videos. For some general principles on the tool, we refer for example to [1].

## 2. Specific usefulness for mathematics

Mathematics is widely recognized as a fundamental domain of knowledge. It is important in most parts of industry and innovation (finance, medical research, technology, earth sciences…), but is also an important part of the development of fundamental skills: abstraction, reasoning or familiarity with important theoretical concepts (numbers, charts, probability…). It is also a domain for which the video technology can be particularly fit since it allows to present a lot of things to schoolchildrens and students in ways that would be extremely difficult to reproduce in a manual or in a classroom. For example, the handling of long algebraic expressions, or the presentation of carefully designed geometrical figures, can be considerably eased by resorting to a video showing these in a dynamical presentation alongside a presentation of their important aspects. Also, preparing the presentation of a relevant and visual example for the solution of some exercise can be a quite long and difficult task for a teacher. In many occasions, this could be advantageously replaced by the use of an already existing video. Videos also allow to present *mathematical experiments* that can be impossible, unconvincing or hazardous to make in the classroom, depending on resources available (paywalled softwares, school materials…).

Hence, even if videos should definitely not be regarded as a kind of universal or magical tool, it can significantly extend the realm of pedagogic opportunities to present mathematics in an enlightening way.

Some authors and teachers argue that the time needed to prepare a video can be minimal, since most of us just do not have enough time to spend to prepare the script of a video (which is, indeed, a long and difficult task), as for example [2]. It is our opinion that such an effort has still to be made since, if successful, will provide a kind of optimal tool that could be useful beyond videomaker's students.

## 3. An academic playlist from Morocco and France for differential equations

In 2021, a project of the authors won a tender of the French Embassy in Morocco and the Moroccan Ministère de l'Éducation Nationale about digital teaching at a university level. The project consisted in a series of 12 short videos (between 3 and 8 minutes each), in French. The subject was ordinary differential equations, which is a standard chapter of most scientific bachelor's degree curricula. These videos are now freely avaliable on-line, from both the websites of the département de mathématiques of the Sorbonne Paris North University [5] and the Moroccan LaR2A (Tétouan Faculty of Sciences) [6] as well as on a dedicated YouTube chain [7]. It is therefore an academic playlist, endorsed by two universities of two different countries. From the beginning, this point led the project to be build up with scientific reliability as its fundamental focus. The authors also tried to be very careful about the two pitfalls previously mentioned (video as a classroom lesson or as a manual) as well as the temptation to turn these videos into some « mathematical entertainment » that would lack the seriousness and/or the scientific content which is to be expected from an institutional project for scientific training.

The institutional aspect of the project also implied that the videos would not be reduced to some catalog of purely technical resolutions of standard exercises. Even if this kind of videos is much appreciated by students, the primary aim was to provide a significant intellectual content: proofs, explanations, links with other sciences, historical aspects, interest for industry, open questions, etc. It could therefore be useful for the training of teachers as well.

Originally, the videos were intended to be a part of a broader project of a website offering written complementary explanations as well as exercises and Moodle-like evaluations. Time as well a financing did not allow all of this to materialize, but such kind of association of different digital tools seems a highly promising possiblility to be set up in the future.

## 4. Some technico-pedagogical aspects

The videos of the project are made of five different kind of shots, shown in figure 1.

The first one, say "institutional", appears only in the beginning of each video, for a period of about 4 seconds. It presents the sponsors of the project and is intended to emphasize on the academic and institutional aspect of the videos.

The second shot, the "human" one, with only the presenters on the screen, is used for introductions, conclusions, transitions, or general recommendations with no mathematical technicalities (e.g. the importance of checking the result of any complex calculation).

The third kind of shot, more technical in its conception, which shows the two presenters together with a small picture somewhere in the screen, is frequently used as a kind of remembering of some result previously given, without explicit statement. The small picture may be assimilated to a comics speech bubble that shows something coming back to the mind of a character.

The fourth kind of shot, centered on essentail content, is used for mathematical technicalities (statement of a theorem, beginning of a calculation…) when they do not (still) fill the whole screen, thus allowing to keep a « human » part on it. One may argue that this could constitue a cognitive load, since it is not strictly necessary to the understanding and also a kind of distraction. Anyway, such a shot is directly followed by another of the last kind, in which presenters do not appear anymore. This let the screen become the equivalent of a blackboard, on which the hard work is done, or some relevant chart is presented, or some visual animation is shown. Whenever possible, this latter shot is animated.

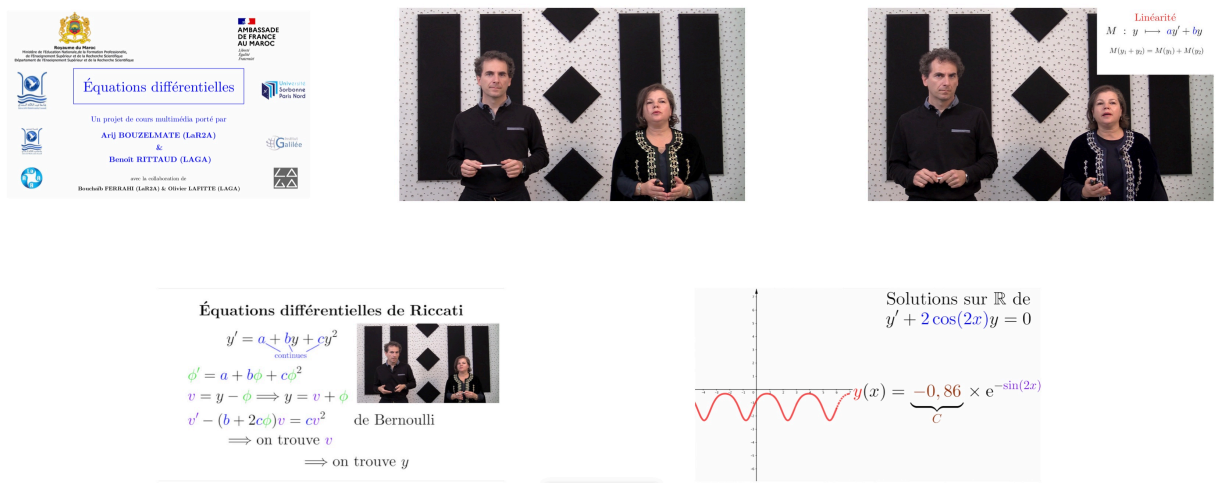

Figure 1: The five kinds of shots: institutional, human, with movie techniques, to the point, and animated.

As already reported, the usefulness of a video for teaching mathematics relies highly on specific opportunities provided by digital technology. In particular, animations are very important (especially when the shot is in « blackboard » mode), hence are used as often as possible. The most entrertaining ones are those showing some visual statements, as in the case of De Beaune problem (find a curve such that, for any of its points, the distance between the abscissa and the intersection of the tangent with the *x* axis is constant). An animation is very enlightening to understand the problem as well as its modelization in analytic terms (figure 2).

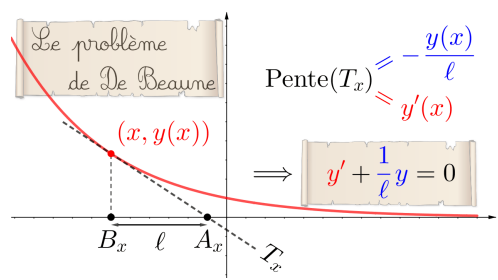

Figure 2: Animation of De Beaune's problem (find a curve such that the distance between the vertical projection of any of its points to the intersection of its tangent with the *x*-axis is constant). In the video, the point $(x, y(x))$ is moving, and the text and relevant lines (some not shown here) appear step by step.

Nevertheless, in general, a mathematical course does not offer that often the possibility of showing something which is entertaining and mathematically substantial at the same time. Therefore, animation refers more often to the successive steps of a reasoning or a calculation, for which progressivity is of utmost importance in a pedagogical perspective. Here are some of the successive screens that appear in the video that presents the implementation, on an explicit example, of a common method to solve some type of linear differential equations (figure 3). Since the aim of the present paper is directly relied to the mathamatical pedagogy, we will not enter here in the details of the substance, restricting ourselves to the form of the slides. (Also, we do not show all the involved screens till the last one, 17 of them being involved in the video.)

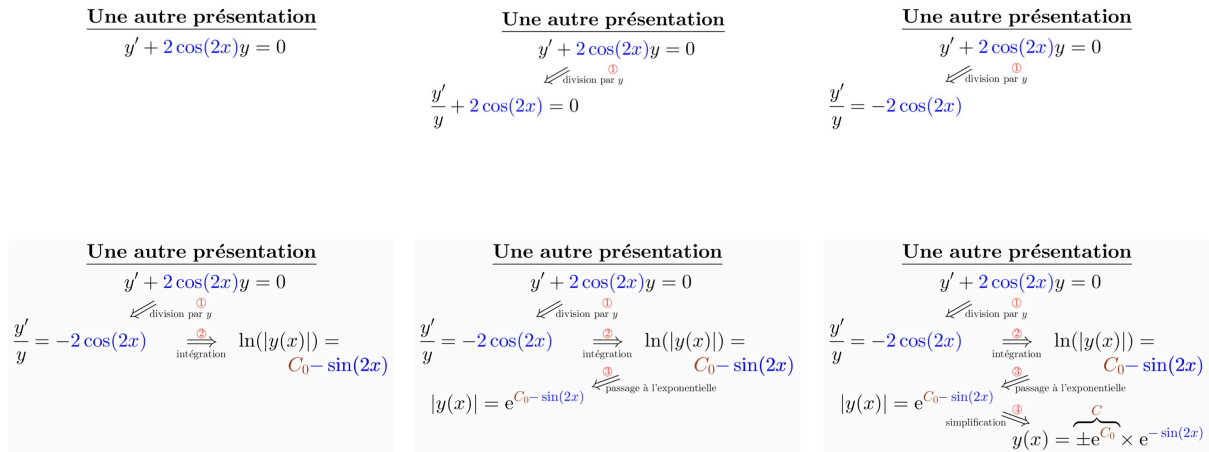

Figure 3: A part of the set of successive screens presenting a way to solve differential equations.

The different elements appear on the screen as the voice-over provides explanation. The differents parts appear step by step to help the understanding since, as also noticed in [3], this way is more efficient than the one in which everything is put on the screen from the beginning. (To be precise, in [3] is compared two kinds of presentations, in which the teacher itself is either writing on the board or elaborating on already drawn illustrations.)

The first most important constraint is that the slides should be readable on the screen of a smartphone, since, as it appeared to many teachers during the pandemics and the lockdowns, many students or schoolchildrens do not have an easy home access to a computer. For pedagogical purposes, it is also reasonable to require that the full resolution of the exercise appears entirely on the final slide. Since the presented method involves several intermediate steps, the screen shows mathematical expressions making zigzags, to make use of all possible space. This could give the impression that the method itself is somehow messy. Hence, to suggest a kind of « backbone » for it, the four arrows indicating the path of resolution are labelled so that the numbers (in red) appear in a single vertical line, at equal distances from each other.

Different colours are used in mathematical expressions, here and in all the videos, in a way as coherent as possible: integration constants $C$ or $C_0$ always appear in brown, terms arising from the coefficient before $y$ (the unknown function to be found) in blue, etc. This was intended to help the watcher to catch more straightforwardly where this or that term comes from. Such a help is especially important since the smallness of the screen did not, in general, allow more than a few explanation words to be written. Also about colours, it is sometimes difficult to decide which one should be used. If a blue 2 is to be added to a red 5, it is quite natural to write the result, 7, in purple. But then, if we add to this 7 some other number written in some other colour, things become more tricky. Studies should be made to provide at least some general and coherent answers to this kind of questions.

More subtle aspects of the previous slides could be presented here, but for the sake of brievety we will mention only one more, about the second and third slides. These ones differ by the single fact that, in the last one, the term 2cos(2$x$) is on the right side of the equality. This is intended to show explicitely the two steps of the calculation (division by $y$, then isolation of $y'/y$), but also to indicate to the student that he has to accept that these two steps are trivial enough to be merged into a single operation henceforth.

Videos also allow to present mathematics in a way that combines written and oral characteristic of a statement. For example, here is how the main theorem for the resolution of a differential equation of the form $y'+by = 0$ is presented, the original oral text being in French (figure 4).

| | |
|---|---|
| The first kind of differential equations we are going to study is the linear differential equations of order one, homogeneous and with coefficient 1 for $y'$, that is: of the form $y'+by = 0$, | 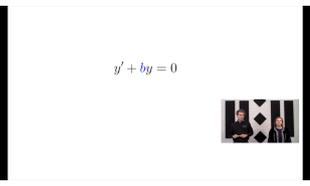 |
| where $b$ is some continuous function defined on an interval $I$. | 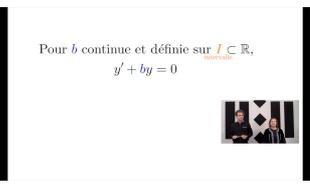 |
| The solutions for this equation are given by a general theorem: | 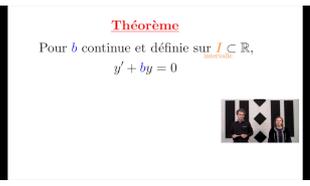 |
| these are the functions $y$ defined on $I$ and of the form $y(x) = C \times \exp(-F(x))$, | 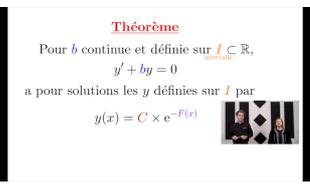 |
| where $C$ is a real constant and $F$ an antiderivative of $b$. | 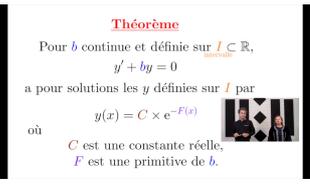 |

Figure 4: The oral sequence of the presentation of a theorem.

As we see, the successive elements of the statement do not appear from left to right and from top to bottom as we could expect. Instead, it is done in a way suited to match better what corresponds to some « natural » oral presentation. Hence, the timeline obeys the rule of such an oral presentation but, in the end, the final statement still appears in a standard written form. This was an attempt to retain the benefits of both oral and written presentation: clarity for the first one and precision for the second.

## 5. Conclusion: new tools, new rules

A lot of work remains to be done in the field of video for mathematical teaching. Some didactical aspects about the ones presented here were discussed in the symposium ETM7 [4], after which it appeared that the understanding of the subject is still in infancy but very promising, for concrete use as well as from a theoretical standpoint. Most importantly, to increase quality as well as efficiency of the tool, forthcoming videos should be designed with greater care for filmmaking aspects. This will probably require a cooperation to be set up with a film director, not only to make the production more professional but also to take better in consideration the filmic aspects that could reinforce the efficiency of the videos.